\documentclass[12pt]{article}
\usepackage{amssymb}
\usepackage{amsfonts}
\usepackage{graphicx}
\usepackage{amsmath}
\usepackage{booktabs}
\usepackage{multirow}
\usepackage{listings}
\usepackage{xcolor}

\definecolor{codegreen}{rgb}{0,0.6,0}
\definecolor{codegray}{rgb}{0.5,0.5,0.5}
\definecolor{codepurple}{rgb}{0,0,1}
\definecolor{backcolour}{rgb}{0.95,0.95,0.92}

\lstdefinestyle{mystyle}{
    backgroundcolor=\color{backcolour},   
    commentstyle=\color{codegreen},
    keywordstyle=\color{blue},
    numberstyle=\tiny\color{codegray},
    stringstyle=\color{black},
    basicstyle=\ttfamily\footnotesize,
    breakatwhitespace=false,         
    breaklines=true,                 
    captionpos=b,                    
    keepspaces=true,                 
    numbers=left,                    
    numbersep=5pt,                  
    showspaces=false,                
    showstringspaces=false,
    showtabs=false,                  
    tabsize=2
}

\newtheorem{theorem}{Theorem}

\newtheorem{algorithm}[theorem]{Algorithm}

\begin{document}
\title{A fast and simple algorithm for the computation of the Lerch transcendent}
\author{Eleonora Denich \thanks{Dipartimento di Matematica e Geoscienze, Universit\`{a} di Trieste, Trieste, Italy, eleonora.denich@phd.units.it} \and Paolo Novati \thanks{Dipartimento di Matematica e Geoscienze, Universit\`{a} di Trieste, Trieste, Italy, novati@units.it}}
\date{}
\maketitle

\begin{abstract}
    This paper deals with the computation of the Lerch transcendent by means of the Gauss-Laguerre formula. An a priori estimate of the quadrature error, that allows to compute the number of quadrature nodes necessary to achieve an arbitrary precision, is derived.
    Exploiting the properties of the Gauss-Laguerre rule and the error estimate, a truncated approach is also considered. 
    The algorithm used and its Matlab implementation are reported.
    The numerical examples confirm the reliability of this approach.
\end{abstract}

\section{Introduction}

In this work we consider the Lerch transcendent introduced in \cite{L} and defined as (see \cite[p.628 n.25.14.1]{NIST})
\begin{equation*}
    \Phi(z,s,a) = \sum_{j=0}^{\infty} \frac{z^j}{(j+a)^s},
\end{equation*}
where $a \neq 0,-1,-2,\ldots$, $|z|<1$ or $\Re(s)>1$, $|z|=1$.
For other values of $z$, $\Phi(z,s,a)$ is defined by analytic continuation.
For a recent investigation on the analytic properties of $\Phi(z,s,a)$, we refer to \cite{LG}. 

The Lerch transcendent and its special cases, such as the polylogarithm (corresponding to $a=1$), appears for instance in quantum Bose and Fermi statistics.
In particular, the particle number density, the pressure, the internal energy and the entropy of the ideal quantum gases of Fermi and Bose can be expressed in terms of the polylogarithm \cite{C}.
Another application of the Lerch transcendent is in biophysics.
Indeed, some discrete distributions, related to the Lerch transcendent by the normalization constant of the associated probability mass functions, are used to establish the statistical composition of DNA and protein sequences \cite{Aksenov}.

In order to evaluate the Lerch transcendent, some authors have used series representations and asymptotic expansions.
Among the others, we quote \cite{Crandall, BB} for an overview of properties, identities and numerical methods for the computation of the Lerch transcendent.
Moreover, in \cite{Aksenov}, combined non-linear condensation transformation, which is an algorithm that allows to evaluate slowly convergent nonalternating series, is considered.
More recently, in \cite{NP}, starting from the Hermite-type integral representation of the Lerch transcendent, a new uniform asymptotic expansion of $\Phi(z,s,a)$, for large order of the parameters $a,s$ and argument $z$, is derived.

In this work, exploiting the integral representation 
\begin{equation} \label{integrale}
\Phi(z,s,a)=\frac{1}{\Gamma(s)}\int_{0}^{\infty}\frac{x^{s-1}e^{-ax}}{1-ze^{-x}}dx, \; z\in \mathbb{C} \setminus [1, + \infty), \Re(s)>0, \Re(a)>0,
\end{equation}
where $\Gamma$ is the Gamma function, we compute $\Phi(z,s,a)$ by employing the Gauss-Laguerre formula and, starting from the analysis given in \cite{B}, we derive an accurate error estimate for the quadrature error.
To the best of our knowledge, the Gauss-Laguerre formula has never been used to this purpose, even if the generalized Laguerre weight clearly appears in the integral representation (\ref{integrale}).
Moreover, after removing the weight, what remains is bounded, and therefore the use of the Gauss-Laguerre method appears a natural choice for this problem.
Finally, since the weights of this rule decay exponentially, we present a reliable algorithm which allows to reduce the number of function evaluations, without significant loss of accuracy in the computation of the integral.

Throughout this work we use the symbol $\approx$ to indicate a generic approximation.
Whereas, the symbol $\sim$ is used to express the asymptotic equality.

The paper is organized as follows.
In Section \ref{section2}, working with $s$ and $a$ real, we present the Gauss-Laguerre approach, together with the error analysis and some numerical experiments.
In Section \ref{section3} we introduce the truncated approach and the algorithm that allows to compute the Lerch transcendent with a prescribed accuracy. 
Some tables with experiments are given in Section \ref{section4}. In Section \ref{section5} we show how to extend the method for $s$ and $a$ complex. An appendix contains an auxiliary result and the Matlab code we have used.

\section{The Gauss-Laguerre approach} \label{section2}

For simplicity we assume $s,a\in \mathbb{R}$. In order to employ the Gauss-Laguerre rule, we consider the change of variable $ax=t$ in (\ref{integrale}), that leads to
\begin{equation*}
\Phi(z,s,a)=\frac{1}{\Gamma(s)a^s}I(z),
\end{equation*}
where
\begin{equation*}
I(z)=\int_{0}^{\infty}\frac{t^{s-1}e^{-t}}{1-ze^{-\frac{t}{a}}}dt.
\end{equation*}
Denoting by 
\begin{equation}\label{f}
I_n(z)=\sum_{j=1}^n w_j^{(n)}f_z\left( t_j^{(n)} \right), \quad f_z(t)=\frac{1}{1-ze^{-\frac{t}{a}}},
\end{equation}
the corresponding $n$-point Gauss-Laguerre rule, where $t_j^{(n)}$ and $w_j^{(n)}$, $j=1,\ldots,n$, are the nodes and the weights respectively, we obtain the approximation
\begin{equation*}
\Phi(z,s,a) \approx \Phi_n(z,s,a)=\frac{1}{\Gamma(s)a^s}I_n(z).
\end{equation*}

\subsection{Error analysis}
In order to derive an estimate for the error 
\begin{equation}\label{En}
E_n= \Phi(z,s,a)-\Phi_n(z,s,a) = \frac{1}{\Gamma(s)a^s} e_n,
\end{equation}
with
\begin{equation}
e_n = I(z)- I_n(z),
\end{equation}
we first need to locate the poles of the function $f_z$ defined in (\ref{f}). By solving $1-ze^{-t/a}=0$, we obtain
\begin{equation}\label{tk}
t_k=a(\ln|z|+i(\arg(z)+2k\pi)),\quad -\pi \leq \arg(z)\leq\pi,\quad k\in\mathbb{Z}.
\end{equation} 
It is known that for the error of the Gauss-Laguerre formula it holds
\begin{equation} \label{espilon}
e_n = \frac{1}{2 \pi i} \int_{\Gamma} \frac{q_n(w)}{L_n^{(s-1)}(w)}f_z(w) dw,
\end{equation}
where $L_n^{(s-1)}$ is $n$-th degree generalized Laguerre polynomial and $q_n$ is the corresponding associated function (see \cite{DR,Szego}). $\Gamma$ is a contour in the complex plane containing the set $[0, + \infty)$ but no singularity of the function $f_z$ can lie on or within $\Gamma$.
For a suitable choice of $\Gamma$, we consider, for $R>1$, the parabola of the complex plane defined by the equation 
\begin{equation}\label{par}
\Re \sqrt{-w} = \ln R.
\end{equation}
It is immediate to verify that the parabola, which we denote by $\Gamma_R$, can be rewritten as 
\begin{equation}\label{para}
w = \frac{y^2}{4 \left( \ln R \right)^2}-\left( \ln R \right) ^2 +i y, \quad y \in \mathbb{R}.
\end{equation}
It is symmetric with respect to the real axis, with vertex in $-\left( \ln R \right) ^2$ and convexity oriented towards the positive real axis. The parabola degenerates to $[0, + \infty)$ as $R\rightarrow1$.

Let $C_k$ be an arbitrary small circle surrounding the pole $t_k$. Assuming $-\pi<\arg(z)<\pi$, the pole closest to the real axis is $t_0$ (see (\ref{tk})), and therefore for any non negative integer $N$ we can define 
\begin{equation} \label{Gamma}
\Gamma=\Gamma_R\bigcup\left(\bigcup_{k=-N}^N C_k\right),
\end{equation}
with $R$ such that the parabola $\Gamma_R$ contains in its interior the poles $t_{-N},...,t_N$ but none of the others.

\begin{figure}
\begin{center}
\includegraphics[scale=0.33]{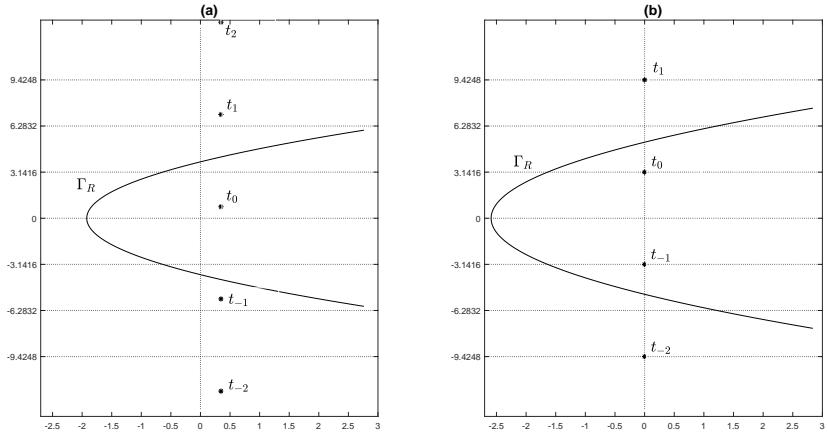}
\end{center}
\caption{The parabola $\Gamma_{R}$ and the poles of the function $f_z(t)$ for $z=\sqrt{2}e^{i \frac{\pi}{4}}$ (left) and $z=e^{i\pi}$ (right). In both cases $a=1$.}
\label{parabole}
\end{figure}

In order to evaluate (\ref{espilon}), with $\Gamma$ as in (\ref{Gamma}), we first observe that, for n $\rightarrow +\infty$, \cite[p.33]{DR}
\begin{equation} \label{q/L}
\frac{q_n(w)}{L_n^{(s-1)}(w)} \sim -2 e^{-i \pi (s-1)} w^{s-1} e^{-w} \frac{K_{s-1}\left(2 \sqrt{m} e^{-i \pi /2} \sqrt{w} \right)}{I_{s-1}\left(2 \sqrt{m} e^{-i \pi /2} \sqrt{w} \right)}, \quad 
\end{equation}
where $m=n+\frac{s}{2}$ and $I$, $K$ are the modified Bessel functions of the first and second kind, respectively. 
Since $0 < \arg(w) < 2 \pi$, we can use for $\alpha>-1$ the asymptotic expansions (see \cite[p.377 n.9.7.1-2]{Abramowitz})
\begin{align*}
I_{\alpha}(t) &= \frac{e^t}{\sqrt{2 \pi t}} \left( 1+ \mathcal{O} \left(\frac{1}{t} \right) \right), \\
K_{\alpha}(t) &= e^{-t} \sqrt{\frac{\pi}{2t}} \left( 1+ \mathcal{O} \left(\frac{1}{t} \right) \right),
\end{align*} 
valid for large $|t|$, $\left\vert \arg (t) \right\vert < \frac{\pi}{2}$, and therefore, 
\begin{equation} \label{K/I}
\frac{K_{\alpha}(t)}{I_{\alpha}(t)} \sim \pi e^{-2t} \left( 1+ \mathcal{O} \left(\frac{1}{t} \right) \right).
\end{equation}
By using (\ref{K/I}) in (\ref{q/L}), we find
\begin{equation}\label{asy}
\frac{q_n(w)}{L_n^{(s-1)}(w)} \sim 2 \pi e^{-i s \pi} w^{s-1} e^{-w} e^{-4 \sqrt{m} \sqrt{-w}}, \quad n \rightarrow + \infty,
\end{equation}
and therefore 
\begin{equation*}
e_n \sim - i e^{-i s \pi} \int_{\Gamma} w^{s-1} e^{-w} e^{-4 \sqrt{m} \sqrt{-w}} f_z(w) dw .
\end{equation*}
By defining for simplicity 
\begin{equation*}
g(w) :=w^{s-1} e^{-w} e^{-4 \sqrt{m} \sqrt{-w}},
\end{equation*}
and with $\Gamma$ as in (\ref{Gamma}), we obtain
\begin{equation}
e_n \sim - i e^{-i s \pi} \left[ \int_{\Gamma_R} g(w) f_z(w) dw -2 \pi i \sum_{k=-N}^N {\rm Res} \left(g(w) f_z(w), t_k \right) \right],
\end{equation}
where ${\rm Res(\cdot,\cdot)}$ denotes the residue. At this point, by (\ref{asy}) and (\ref{par})
\begin{eqnarray*}
\int_{\Gamma_R} g(w) f_z(w) dw &=& \int_{\Gamma_R} w^{s-1} e^{-w} e^{-4 \sqrt{m} \sqrt{-w}} f_z(w) dw \\
&=& R^{-4 \sqrt{m}} \int_{\Gamma_R} w^{s-1} e^{-w} \chi(w) f_z(w) dw,
\end{eqnarray*}
where 
\begin{equation} \label{chi}
 \chi(w)=e^{-4 i\sqrt{m} \Im\sqrt{-w}}.   
\end{equation}
We remark that, since $|f_z|$ is bounded, the above integral is also bounded.
Moreover, since  ${\rm Res} \left(f_z(w), t_k \right)=a$, $k\in\mathbb{Z}$, by (\ref{tk}) we have
\begin{eqnarray*}
{\rm Res} \left(g(w)f_z(w), t_k \right) &=& g(t_k){\rm Res} \left(f_z(w), t_k \right) \\
&=& a^s \left[{\rm ln}_k(z)\right]^{s-1} z^{-a}\chi(t_k) R_k^{-4\sqrt{m}},
\end{eqnarray*}
where we use the notation
\begin{equation*}
{\rm ln}_k(z)=\ln|z|+i(\arg(z)+2k\pi), \quad k\in\mathbb{Z},
\end{equation*}
and $R_k<R$ (for $-N\leq k \leq N$) represents the parabola passing through $t_k$, that is, $R_k$ solves $\Re \sqrt{-t_k} = \ln R_k $. 
Joining the above results we finally obtain
\begin{equation*} 
\begin{split}
e_n &\sim - i e^{-i s \pi} \left[R^{-4 \sqrt{m}} \int_{\Gamma_R} w^{s-1} e^{-w} \chi(w) f_z(w) dw \right. \\
&-2\pi i a^s z^{-a} \left. \sum_{k=-N}^N   \left[{\rm ln}_k(z)\right]^{s-1} \chi(t_k) R_k^{-4\sqrt{m}} \right].
\end{split}
\end{equation*}
Under the assumption $\left\vert \arg (z) \right\vert < \pi$, we have $R_0<R_k<R$, for $|k|=1,...,N$, and therefore, by collecting the term involving $R_0$,
\begin{equation} \label{e_n_finale}
e_n \sim -2\pi e^{-i s \pi}  a^s \left[{\rm ln}_0(z)\right]^{s-1} z^{-a}\chi(t_0) R_0^{-4\sqrt{m}}.
\end{equation}
This situation in shown in Figure \ref{parabole}a, where $z=\sqrt{2}e^{i \frac{\pi}{4}}$.
In the case of $\arg(z)= \pm \pi$ we have that the poles are symmetric with respect to the real axes and it holds $\overline{t_{k}}=t_{-k-1}$ if $\arg(z)=\pi$, whereas $\overline{t_{-k}}=t_{k+1}$ if $\arg(z)=-\pi$, $k \geq 0$ (see Figure \ref{parabole}b, in which $z=e^{i\pi}$).
In both situations, formula (\ref{e_n_finale}) needs to be replaced by the sum of the terms corresponding to the poles $t_0$ and $t_{-1}$, for $\arg(z) = \pi$, and $t_0$ and $t_1$ for $\arg(z) =-\pi$. In both cases, we obtain the same result. 
Assuming for instance $\arg(z)=\pi$, we have that $\overline{t_0}=t_{-1}$, $R_0=R_{-1}$ and therefore
\begin{eqnarray*}
e_n &\sim& -2\pi e^{-i s \pi}  a^s z^{-a} R_0^{-4\sqrt{m}} \left[\left[{\rm ln}_0(z)\right]^{s-1} \chi(t_0) + \left[{\rm ln}_{-1}(z)\right]^{s-1} \chi(t_{-1}) \right]\\
&=& -4\pi e^{-i s \pi}  a^s z^{-a} R_0^{-4\sqrt{m}} \Re \left[\left[{\rm ln}_0(z)\right]^{s-1} \chi(t_0)  \right]. 
\end{eqnarray*}
Working with the modulus, the only essential difference between the above formula and (\ref{e_n_finale}) is a factor 2. Moreover, since the latter should also be used for $\arg(z)$ close to $\pm\pi$, independently of $z$ we finally consider the error estimate
\begin{align} \label{epsilon_n_finale}
|e_n| &\approx 4\pi a^s |z|^{-a}  |{\rm ln}_0(z)|^{s-1} R_0^{-4\sqrt{m}} \notag \\
& =4\pi a^s |z|^{-a}  |{\rm ln}_0(z)|^{s-1} e^{-4 \sqrt{m}\Re \left( \sqrt{-t_0}\right) } =: \epsilon_n,
\end{align}
for the integral, and therefore (see (\ref{En}))
\begin{equation}\label{error}
|E_n| \approx \frac{4\pi}{\Gamma(s)} |z|^{-a}  |{\rm ln}_0(z)|^{s-1} e^{-4 \sqrt{m}\Re \left( \sqrt{-t_0}\right) } =: \mathcal{E}_n.
\end{equation}

\subsection{Some numerical experiments} \label{section2.2}

\begin{figure}
\begin{center}
\includegraphics[scale=0.33]{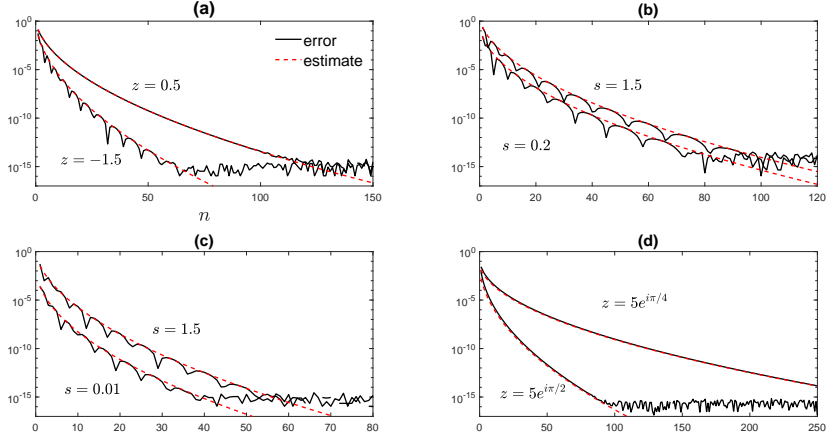}
\end{center}
\caption{Error and error estimate (\ref{error}) in the case of (a) the polylogarithm for $s=3/2$, (b) the Dirichlet beta function, (c) the Dirichlet eta function and (d) the Lerch transcendent for $s=1.5$ and $a=2.5$.}
\label{figure_ee}
\end{figure}

In Figure \ref{figure_ee} we show the quality of estimate (\ref{error}) on some examples.
Specifically, we consider the polylogarithm \cite[n.25.12]{NIST} (Figure \ref{figure_ee}a)
\begin{equation*}
    \text{Li}_s(z) = \sum_{j=1}^{\infty} \frac{z^j}{j^s},
\end{equation*}
the Dirichlet beta function \cite[p.807, n.23.2.21]{Abramowitz} (Figure \ref{figure_ee}b)
\begin{equation*}
    \beta(s)=\sum_{j=0}^{\infty} \frac{(-1)^j}{(2j+1)^s},
\end{equation*}
and the Dirichlet eta function \cite[p.807, n.23.2.19]{Abramowitz} (Figure \ref{figure_ee}c)
\begin{equation*}
    \eta(s)=\sum_{j=1}^{\infty} \frac{(-1)^{j-1}}{j^s}.
\end{equation*}
These functions can be written in terms of $\Phi(z,s,a)$ for particular values of the parameters, that is,
\begin{align*}
    \text{Li}_s(z) &= z \Phi(z,s,1), \\
    \beta(s) &= 2^{-s} \Phi(-1,s,1/2), \\
    \eta(s) &= \Phi(-1,s,1).
\end{align*}
We remark that the value $s=3/2$ chosen for the polylogarithm represents the case of the three-dimensional Bose gas in a box.
Finally, in Figure \ref{figure_ee}d we consider a general situation.
As for the computation of the nodes and weights of the generalized Gauss-Laguerre rule (c.f. (\ref{f})), we employ the Matlab routine $\mathtt{lagpts.m}$ \cite{T}, based on the traditional Golub-Welsch algorithm (see \cite{GW,DR}).
We point out that the oscillations of the error in Figure \ref{figure_ee}a-b-c are due to the term $\chi(t_0)$ (cf. (\ref{chi})- (\ref{e_n_finale})).

\section{A truncated approach} \label{section3}

Having at disposal an accurate error estimate, in this section we derive a reliable algorithm that allows to reduce the number of evaluations of the function $f_z$ (see (\ref{f})), required by the quadrature rule. This is possible because the weights of the Gauss-Laguerre rule decay exponentially and, moreover, because the function $f_z$ is bounded. 
Indeed, a simple analysis, reported in Appendix, shows that
\begin{equation}\label{Kz}
|f_z(t)| = \left\vert 1-ze^{-t/a} \right\vert^{-1} \leq K_z:=
\begin{cases}
1& \quad \Re (z) \leq 0 \\
\frac{|z|}{\left\vert \Im (z) \right\vert}& \quad \Re (z) > 0, \Re(z) \leq |z|^2 \\
|1-z|^{-1}& \quad \Re(z)>0, \Re(z) > |z|^2
\end{cases}.
\end{equation}
Now, let $\tau_n$ be the solution of 
\begin{equation}\label{liminf}
K_z \int_{\tau_n}^{+\infty} t^{s-1}e^{-t} dt = \epsilon_n,
\end{equation}
where $\epsilon_n$ is defined in (\ref{epsilon_n_finale}). By using the relation \cite[p.942 n.8.357]{GR}
\begin{equation*}
\int_{\tau_n}^{+\infty} t^{s-1}e^{-t} dt = \tau_n^{s-1} e^{-\tau_n} \left( 1+ \mathcal{O} \left( \frac{1}{\tau_n}\right)\right),
\end{equation*}
we approximate $\tau_n$ by solving with respect to $x$ the equation
\begin{equation*}
x^{s-1} e^{-x} = \frac{\epsilon_n}{K_z}.
\end{equation*}
The solution is given by
\begin{equation*}
x=(1-s) W\left( \frac{1}{1-s} \left( \frac{\epsilon_n}{K_z} \right)^{\frac{1}{s-1}} \right),
\end{equation*}
where $W$ denotes the Lambert W function. In particular, we obtain 
\begin{equation}
x=(1-s) W_0\left( \frac{1}{1-s} \left( \frac{\epsilon_n}{K_z} \right)^{\frac{1}{s-1}} \right), 
\end{equation}
for $0<s<1$, and
\begin{equation}
x=(1-s) W_{-1}\left( \frac{1}{1-s} \left( \frac{\epsilon_n}{K_z} \right)^{\frac{1}{s-1}} \right),
\end{equation}
for $s>1$.
It is known that (\cite[n. 4.13.10-11]{NIST})
\begin{align*}
W_0(y) &\sim \ln(y), \quad {\rm for} \quad y \rightarrow + \infty, \\
W_{-1}(-y) &\sim -\ln(y), \quad {\rm for} \quad y \rightarrow 0^+,
\end{align*}
and therefore $x \sim g_n(s)$, where
\begin{equation} \label{prima}
g_n(s) = - \ln\left( \frac{\epsilon_n}{K_z} \right) +(s-1) \ln |1-s|.
\end{equation}
Note that $g_n(s)$ has a removable singularity in $s=1$, so that we can define
\begin{equation*}
g_n(1) = - \ln\left( \frac{\epsilon_n}{K_z}\right),
\end{equation*} 
that is also the exact solution of (\ref{liminf}), with respect to $\tau_n$, for $s=1$.
In conclusion, we have that $\tau_n \sim g_n(s)$, for $s>0$.

We are now on the point to introduce the truncated rule 
\begin{equation}\label{Ikn}
I_{k_n}(z) = \sum_{j=1}^{k_n} w_j^{(n)} f_z\left( x_j^{(n)} \right) ,
\end{equation}
where $k_n$ is the smallest integer such that $x_j^{(n)} \geq g_n(s)$, for each $j \geq k_n$.
As for the error, we have
\begin{align*}
\left\vert I(z) -I_{k_n}(z) \right\vert &= \left\vert I(z) -I_n(z)+\sum_{j=k_n+1}^n w_j^{(n)} f_z \left( x_j^{(n)} \right) \right\vert \\
& \lesssim \epsilon_n + \sum_{j=k_n+1}^n w_j^{(n)} \left\vert f_z \left( x_j^{(n)} \right) \right\vert \\
&\leq \epsilon_n + K_z \sum_{j=k_n+1}^n w_j^{(n)}.
\end{align*}
Now, using the bound \cite[eqs. 2.4-2.7]{MO}
\begin{equation*}
w_j^{(n)} \leq c \left( x_j^{(n)} - x_{j-1}^{(n)} \right) \left(x_j^{(n)} \right)^{s-1} e^{-x_j^{(n)}},
\end{equation*}
where $c$ is a constant independent of $j,n$ and close to $1$ for large $n$, we have
\begin{align*}
\sum_{j=k_n+1}^n w_j^{(n)} &\leq c \sum_{j=k_n+1}^n \left( x_j^{(n)} - x_{j-1}^{(n)} \right) \left(x_j^{(n)} \right)^{s-1} e^{-x_j^{(n)}} \\
&\leq c \int_{x_{k_n}^{(n)}}^{+\infty} x^{s-1} e^{-x} dx \quad \text{(for $k_n$ large enough)} \\
&\leq c \int_{g_n(s)}^{+\infty} x^{s-1} e^{-x} dx.
\end{align*}
Finally, for the truncated rule we obtain the estimate
\begin{equation} \label{23bis}
\left\vert I(z) -I_{k_n}(z) \right\vert \lesssim (1+c)  \epsilon_n \approx 2 \epsilon_n.
\end{equation}

\begin{figure}
\begin{center}
\includegraphics[scale=0.33]{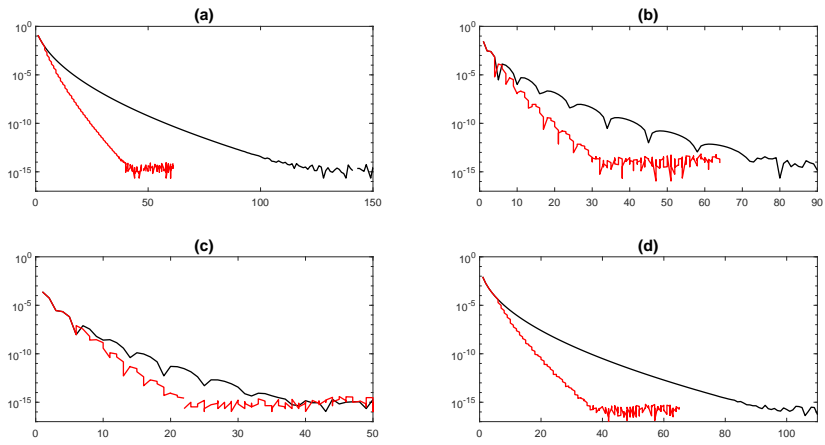}
\end{center}
\caption{Comparison between the error of the generalized Laguerre formula (black) and the truncated approach (red) for (a) the polylogarithm with $s=3/2$ and $z=0.5$, (b) the Dirichlet beta function with $s=0.2$, (c) the Dirichlet eta function with $s=0.01$ and (d) the Lerch transcendent with $s=1.5$, $a=2.5$ and $z=5 e^{i \pi /2}$.}
\label{figure_T}
\end{figure}

In order to understand the decay rate with respect to $k_n$ we use the relations \cite[n. 22.16.8]{Abramowitz} 
\begin{equation*}
x_k^{(n)}=\frac{j_{s-1,k}^2}{4n+2s} \left( 1+ \mathcal{O}\left( \frac{1}{n^2} \right) \right),
\end{equation*}
and \cite[n. 9.5.12]{Abramowitz} 
\begin{equation*}
j_{s-1,k}= \left(k+ \frac{s}{2} -\frac{3}{4} \right) \pi \left( 1+ \mathcal{O}\left( \frac{1}{k} \right) \right),
\end{equation*}
to finally obtain
\begin{equation} \label{zeros}
x_k^{(n)} \approx \frac{k^2 \pi^2}{4m}.
\end{equation}
At this point, it is possible to derive an analytical approximation of $k_n$ by imposing $x_k^{(n)}=g_n(s)$ with the help of (\ref{zeros}). 
In order to detect the new asymptotic behavior, we further simplify the computation by neglecting the term involving $s$ in (\ref{prima}), that is, we solve with respect to $k$
\begin{equation*}
\frac{k^2 \pi^2}{4m} = - \ln\left( \frac{\epsilon_n}{K_z}\right).
\end{equation*}
We remark, however, that this simplification is not used in the algorithm presented below.
By (\ref{epsilon_n_finale}) and defining
\begin{equation*}
C := 4\pi a^s |z|^{-a}  |{\rm ln}_0(z)|^{s-1},
\end{equation*}
we find
\begin{equation*}
 \frac{k^2 \pi^2}{4m} = \ln \left( \frac{K_z}{C} \right) + 4 \sqrt{m} \ln(R_0),
\end{equation*}
so that, for $n$ large enough,
\begin{equation}\label{kn}
k_n \approx  \frac{4}{\pi} m^ \frac{3}{4} \left(\ln(R_0) \right)^\frac{1}{2}.
\end{equation}
By using this relation in (\ref{epsilon_n_finale}) we finally obtain
\begin{equation*}
\epsilon_n \approx C e^{-d\, k_n^\frac{2}{3}}, \quad d=\left(2 \pi \ln(R_0) \right)^\frac{2}{3},
\end{equation*}
which expresses the speed up attainable with the truncation.
Indeed, the number of function evaluations is now raised to the power of $2/3$, whereas in (\ref{epsilon_n_finale}) the power is $1/2$.
This decay rate is well visible in Figure \ref{figure_T}, where we show the benefits of the truncated rule on some examples. In particular, as in Section \ref{section2.2}, we consider the polylogarithm (Figure \ref{figure_T}a), the Dirichlet beta function (Figure \ref{figure_T}b), the Dirichlet eta function (Figure \ref{figure_T}c) and a general situation (Figure \ref{figure_T}d).  

We summarize the basic steps for the computation of the Lerch transcendent by using the truncated Laguerre rule and with a prescribed error tolerance $\mathcal{E}$ in the following algorithm. 

\begin{algorithm}[Truncated Laguerre rule] \label{algoritmo}
Given $z,s,a,\mathcal{E}$
\begin{enumerate}
\item compute the corresponding tolerance for the integral $\epsilon =  a^s \Gamma(s) \frac{\mathcal{E}}{2} $ (see (\ref{En}) and (\ref{23bis}))
\item compute $K_z$ as in (\ref{Kz})
\item compute $R_0$ such that $t_0 \in \Gamma_{R_0}$ by using (\ref{para})
\item solve $\epsilon_n = \epsilon$ with respect to $m$ (see (\ref{epsilon_n_finale})) and then set $n=m-s/2$
\item estimate $k_n$ by solving, with respect to $k$, $x_k^{(n)} = g_n(s)$, by using (\ref{zeros}) and (\ref{prima}) (in the implementation it is convenient to add one or two more points in order to prevent the under estimation of the error, because of the large number of approximations used)
\item compute the first $k_n$ nodes and weights $x_j^{(n)},w_j^{(n)}$, $j=1,...,k_n$, of the $n$-point Laguerre rule
\item approximate $\Phi(z,s,a)$ with $\Phi_{k_n}(z,s,a)=\frac{1}{\Gamma(s)a^s} I_{k_n}(z)$ (see (\ref{Ikn}))
\end{enumerate}
\end{algorithm}

\section{Numerical experiments} \label{section4}

In this section, working with prescribed error tolerances $\mathcal{E}=1e-10$ and $\mathcal{E}=1e-14$, we test Algorithm \ref{algoritmo} on several examples.
In particular, we present some tables in which, for each set of parameters, we report the corresponding values of $n$ and $k_n$, given by steps $4$-$5$, together with the final error obtained by using a reference solution.
Specifically, we consider the polylogarithm $\text{Li}_s\left( r e^{i \tau \pi} \right)$ in Table \ref{table1}, the Dirichlet beta function $\beta(s)$ in Table \ref{table2}, the Dirichlet eta function $\eta(s)$ in Table \ref{table3} and general cases of the Lerch transcendent $\Phi\left( r e^{i \tau \pi},s,a \right)$ in Table \ref{table4}, for different values of the parameters $r, \tau,s,a$.
In Tables \ref{table1} and \ref{table4}, $\tau$ is set in order to consider the arguments $\pi, \frac{3}{4} \pi, \frac{\pi}{2}, \frac{\pi}{4}$.
We can see that, except in few rare cases, the prescribed tolerance is achieved.

\begin{table} 
\[
\begin{array}{ccccccccc}
\toprule
  &      &    & \multicolumn{3}{c}{\mathcal{E}=1e-10} & \multicolumn{3}{c}{\mathcal{E}=1e-14} \\
r & \tau & s  & n & k_n & {\rm error}  & n & k_n & {\rm error}\\
\midrule
0.5 & 1 & 1.5  & 25 & 18 & 2.49e-11 & 44 & 27 & 1.06e-15 \\
 & 0.75 &       & 31 & 20 & 1.21e-11 & 54 & 29 & 2.11e-15 \\
  & 0.5 &       & 39 & 22 & 1.78e-11 & 70 & 33 & 3.83e-15 \\
  & 0.25 &       & 53 & 25 & 1.89e-11 & 95 & 38 & 1.62e-15 \\ 
  \midrule
2 & 1 & 1.5  & 35 & 21 & 1.13e-11 & 63 & 31 & 4.22e-15 \\  
 & 0.75 &      & 49 & 24 & 1.90e-11 & 89 & 37 & 4.46e-15 \\
  & 0.5 &       & 83 & 31 & 2.09e-11 & 151 & 47 & 3.98e-15 \\
  & 0.25 &       & 233 & 50 & 2.58e-11 & 430 & 78 & 4.91e-15 \\ 
  \midrule
0.7 & 1 & 0.5  & 24 & 18 & 1.73e-11 & 43 & 26 & 7.66e-15 \\
   & 0.75 &       & 30 & 19 & 2.28e-11 & 56 & 30 & 1.51e-15 \\
  & 0.5 &       & 43 & 23 & 1.90e-11 & 78 & 35 & 3.24e-15 \\
  & 0.25 &       & 70 & 29 & 2.34e-11 & 128 & 44 & 6.24e-15 \\
  \midrule
3 & 1 & 0.5  & 33 & 20 & 2.95e-11 & 62 & 31 & 8.27e-15 \\
   & 0.75 &       & 49 & 24 & 2.22e-11 & 93 & 38 & 1.51e-15 \\
  & 0.5 &       & 90 & 32 & 2.48e-11 & 172 & 50 & 2.33e-14 \\
  & 0.25 &       & 297 & 56 & 2.09e-11 & 562 & 89 & 2.51e-15 \\ 
\bottomrule
\end{array}
\]
\caption{Results of Algorithm \ref{algoritmo} for the polylogarithm $\text{Li}_s\left( r e^{i \tau \pi} \right)$.} \label{table1}
\end{table}

\begin{table} 
\[
\begin{array}{ccccccc}
\toprule
   & \multicolumn{3}{c}{\mathcal{E}=1e-10} & \multicolumn{3}{c}{\mathcal{E}=1e-14} \\
s  & n & k_n & {\rm error} & n & k_n & {\rm error}\\
\midrule
0.5 & 51 & 25 & 3.44e-11 & 94 & 38 & 1.99e-14 \\
1 & 55 & 26 & 4.20e-11 & 101 & 40 & 7.55e-15 \\
1.5 & 58 & 27 & 3.94e-11 & 105 & 40 & 5.33e-15\\
2 & 59 & 27 & 3.65e-11 & 106 & 41 & 1.33e-15 \\
2.5 & 60 & 28 & 1.84e-11 & 108 & 42 & 8.88e-16 \\
3 & 61 & 28 & 2.56e-12 & 109 & 42 & 4.44e-15 \\
3.5 & 61 & 29 & 3.57e-11 & 109 & 43 & 1.78e-15 \\
4 & 60 & 29 & 6.83e-11 & 108 & 43 & 5.33e-15 \\
4.5 & 60 & 29 & 5.18e-11 & 108 & 44 & 3.55e-15 \\
5 & 59 & 30 & 9.07e-12 & 107 & 44 & 7.11e-15 \\
\bottomrule
\end{array}
\]
\caption{Results of Algorithm \ref{algoritmo} for the Dirichlet beta function $\beta(s)$.} \label{table2}
\end{table}

\begin{table} 
\[
\begin{array}{ccccccc}
\toprule
& \multicolumn{3}{c}{\mathcal{E}=1e-10} & \multicolumn{3}{c}{\mathcal{E}=1e-14} \\
s  & n & k_n & {\rm error} & n & k_n & {\rm error}\\
\midrule
0.5 & 26 & 18 & 3.28e-12 & 47 & 27 & 1.37e-14 \\
1 & 28 & 19 & 5.82e-12 & 51 & 29 & 2.22e-16 \\
1.5 & 29 & 19 & 8.89e-13 & 53 & 29 & 1.22e-15 \\
2 & 29 & 19 & 3.36e-11 & 53 & 29 & 2.66e-15 \\
2.5 & 30 & 20 & 3.21e-11 & 54 & 30 & 4.33e-15 \\
3 & 30 & 20 & 4.78e-11 & 54 & 30 & 4.22e-15 \\
3.5 & 30 & 20 & 4.39e-11 & 54 & 30 & 3.44e-15 \\
4 & 29 & 20 & 5.02e-12 & 53 & 30 & 1.22e-15 \\
4.5 & 29 & 20 & 3.29e-11 & 53 & 30 & 4.11e-15 \\
5 & 29 & 21 & 5.25e-11 & 53 & 31 & 5.33e-15 \\
\bottomrule
\end{array}
\]
\caption{Results of Algorithm \ref{algoritmo} for the Dirichlet eta function $\eta(s)$.} \label{table3}
\end{table}

\begin{table} 
\[
\begin{array}{cccccccccc}
\toprule
  &      &   &    & \multicolumn{3}{c}{\mathcal{E}=1e-10} & \multicolumn{3}{c}{\mathcal{E}=1e-14} \\
r & \tau & s & a  & n & k_n & {\rm error} &  n & k_n & {\rm error}\\
\midrule
0.5 & 1 & 0.5 & 0.7 & 30 & 19 & 4.11e-11 & 56 & 30 & 3.33e-16\\
 & 0.75 &     &  & 38 & 22 & 1.97e-11 & 70 & 33 & 1.08e-15 \\
  & 0.5 &     &  & 50 & 24 & 2.16e-11 & 92 & 38 & 6.97e-15 \\
  & 0.25 &    &   & 71 & 29 & 2.01e-11 & 129 & 44 & 7.65e-15 \\ 
\midrule
2 & 1 & 1.4 & 2 & 17 & 15 & 3.19e-12 & 30 & 22 & 3.91e-15 \\
 & 0.75 &     &  & 23 & 17 & 2.12e-11 & 43 & 26 & 1.48e-15 \\
  & 0.5 &     &  & 39 & 21 & 2.17e-11 & 73 & 33 & 2.71e-15 \\
  & 0.25 &    &   & 111 & 35 & 2.14e-11 & 207 & 54 & 5.27e-15 \\ 
\midrule
5 & 1 & 0.2 & 1.1 & 29 & 19 & 5.07e-11 & 58 & 30 & 5.33e-15 \\
 & 0.75 &     &  & 45 & 23 & 2.61e-11 & 90 & 36 & 1.46e-14 \\
  & 0.5 &     &  & 89 & 31 & 2.45e-11 & 177 & 50 & 3.11e-14 \\
  & 0.25 &    &   & 317 & 57 & 2.35e-11 & 630 & 93 & 2.17e-15 \\ 
\midrule
8 & 1 & 4 & 3 & 11 & 11 & 4.88e-11 & 23 & 20 & 1.16e-14 \\
 & 0.75 &     &  & 17 & 15 & 5.88e-11 & 36 & 24 & 4.46e-15 \\
  & 0.5 &     &  & 34 & 20 & 3.56e-11 & 71 & 32 & 3.18e-15 \\
  & 0.25 &    &   & 122 & 35 & 5.70e-11 & 257 & 59 & 3.37e-14 \\ 
\bottomrule
\end{array}
\]
\caption{Results of Algorithm \ref{algoritmo} for general cases of the Lerch transcendent $\Phi\left( r e^{i \tau \pi},s,a \right)$.} \label{table4}
\end{table}

\section{Complex case} \label{section5}

\begin{figure}
\begin{center}
\includegraphics[scale=0.33]{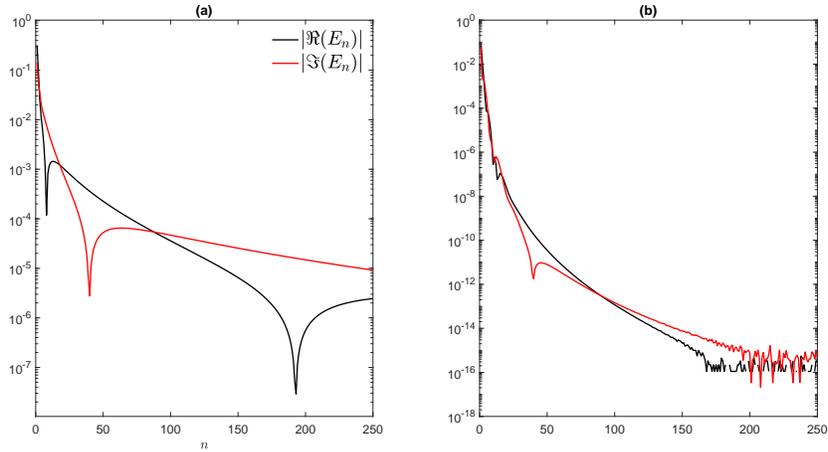}
\end{center}
\caption{Real and imaginary part of the error for $s=2+i$, $a=1$, $z=-1.1$ (left) and $s=8+i$, $a=1$, $z=-1.1$ (right).}
\label{figure_complex}
\end{figure}

Starting from the integral representation (\ref{integrale})
\begin{equation*}
\Phi(z,s,a)=\frac{1}{\Gamma(s)}\int_{0}^{\infty}\frac{x^{s-1}e^{-ax}}{1-ze^{-x}}dx, \; z\in \mathbb{C} \setminus [1, + \infty), \Re(s)>0, \Re(a)>0,
\end{equation*}
by using the change of variable $x=\frac{t}{\Re(a)}$, we obtain 
\begin{equation*}
\Phi(z,s,a)=\frac{e^{-i \Im(s) \ln \Re(a)}}{\Re(a)^{\Re(s)}\Gamma(s)}\int_{0}^{\infty}\frac{t^{\Re(s)-1}e^{-t}}{1-z e^{-i\frac{\Im(a)}{\Re(a)}t}} e^{i\left(\Im(s) \ln t- \frac{\Im(a)}{\Re(a)}t\right)} dt. 
\end{equation*}
The function 
\begin{equation*}
    \varphi(t) = e^{i\left(\Im(s) \ln t- \frac{\Im(a)}{\Re(a)}t\right)}
\end{equation*}
inside the integral is the main difference with respect to the case of $s,a \in \mathbb{R}$.
Some experiments have revealed that the case of $\Im(a) \neq 0$ does not constitute a problem for the Laguerre rule, unless $\Im(a) \gg \Re(a)$.
On the other side, the case of $\Im(s) \neq 0$ may be difficult to handle.
Indeed, as $t \rightarrow 0$, the real and the imaginary part of the function $\varphi(t)$ oscillate with increasing frequency and the Laguerre rule appears to be inadequate (Figure \ref{figure_complex}a).
The situation is overtaken whenever $\Re(s) \gg \left\vert \Im(s) \right\vert$ and $\Re(s)>1$, since the term $t^{\Re(s)-1}$ makes the oscillations negligible near zero (see Figure \ref{figure_complex}b).

\section{Conclusion}

In this work we have employed the generalized Gauss-Laguerre formula to compute the Lerch transcendent $\Phi(z,s,a)$ for $z\in \mathbb{C} \setminus [1, + \infty)$ and $ s,a>0$.
We have derived sharp error estimates that enable to know a priori the number of quadrature points necessary to achieve a prescribed accuracy and to truncate the rule.
We have tested the arising algorithm on several examples and the results confirm the reliability of this approach.
The extension to $s$ and $a$ complex is not theoretically analyzed.
Anyway, at a first glance the Laguerre rule appears robust if the imaginary parts are relatively small with respect to the moduli.
On the contrary, the high frequency of the oscillations makes the method unsuited.

\appendix

\section{Bounds for the integrand function}

In order to derive the bound (\ref{Kz}), we define $x=e^{-t/a}$, $t \in [0,+\infty)$, so that
\begin{equation*}
    |f_z(t)| =|1-zx|^{-1} = \left( 1+|z|^2x^2-2 \Re(z) x \right)^{-1/2}, \quad x \in [0,1].
\end{equation*}
The function $g(x)= 1+|z|^2x^2-2 \Re(z) x$ is a parabola and its minimum in $[0,1]$ defines $K_z$.
The vertex is at $x_V=\frac{\Re(z)}{|z|^2}$ and we have three cases. 
Let $D= \mathbb{C} \setminus [0, + \infty)$ be the domain of definition of integral (\ref{integrale}) with respect to $z$. 
Let $A= \left\lbrace z \in \mathbb{C} \mid \Re(z) \leq 0 \right\rbrace $, $B= \left\lbrace z \in \mathbb{C} \mid \left\vert z-\frac{1}{2}\right\vert <\frac{1}{2} \right\rbrace$ and $C= D \setminus (A \cup B)$.
Note that $A,B,C$ are mutually disjoint and $A \cup B \cup C =D$.
Moreover, the set $B$ is the solution of the inequality $\Re(z)>|z|^2$.
At this point, if $z \in A$, then $x_V \leq 0$ and $g(x)\geq g(0)$, and therefore
\begin{equation*}
    |f_z(t)| \leq 1.
\end{equation*}
If $z \in B$, then $x_V>1$, $g(x) \geq g(1)$ and we find 
\begin{equation*}
    |f_z(t)| \leq |1-z|^{-1}.
\end{equation*}
Finally, if $z \in C$, then $0<x_V \leq 1$, $g(x) \geq g(x_V)$ and we obtain
\begin{equation*} 
    |f_z(t)| \leq \frac{|z|}{|\Im(z)|}.
\end{equation*}

\section{Matlab code}

\begin{lstlisting}[language=Matlab]
function [lerch,kn] = LerchT(z,s,a,Epsilon)
% LERCHT evaluates the Lerch transcendent \Psi(z,s,a) (s,a real), with accuracy
% Epsilon, using the Gauss-Laguerre rule. The result is contained in the 
% output lerch, whereas kn is the number of function evaluations. 
% The code implements the algorithm (Truncated Laguerre rule) of Section 3.

% corresponding tolerance for the integral (step 1)
tol = a^s*gamma(s)*Epsilon/2; 

% computation of K_z (step 2)
if real(z) <= 0
    Kz = 1;
else
    if real(z) <= abs(z)^2
        Kz = abs(z)/imag(z);
    else
        Kz = 1/abs(1-z);
    end
end

% steps 3-4
t0 = a*(log(abs(z))+1i*angle(z)); 
C = 4*pi*a^s*abs(z)^(-a)*abs(t0/a)^(s-1);
m = ceil(((log(C)-log(tol))/(4*real(sqrt(-t0))))^2);
n = ceil(m-s/2);

% step 5
if s == 1
    gn = -log(tol/Kz);
else
    gn = -log(tol/Kz)+(s-1)*log(abs(1-s)); 
end
kn = min(ceil(sqrt(4*m*gn)/pi)+2,n); % number of functions evaluation; the +2 is to avoid under estimation

% step 6
[x,w] = lagpts(n,s-1);
xT = x(1:kn); wT = w(1:kn);

% step 7
A = 1/(gamma(s)*a^s);
f = @(t) (1-z*exp(-t/a)).^(-1);
lerch = A*(wT*f(xT)); 

end
\end{lstlisting}

\section*{Acknowledgements}

This work was partially supported by GNCS-INdAM, FRA-University of Trieste and CINECA under HPC-TRES program award number 2019-04. The authors are member of the INdAM research group GNCS.

\end{document}